\documentclass{amsart}
\usepackage[latin1]{inputenc}
\usepackage{amsmath}
\usepackage{amsthm}
\usepackage{amsfonts}
\usepackage{amssymb}
\newtheorem{theorem}{Theorem}[section]

\newtheorem{proposition}[theorem]{Proposition}
\newtheorem{definition}[theorem]{Definition}

\newtheorem{remark}[theorem]{Remark}
\newtheorem{corollary}[theorem]{Corollary}
\newtheorem{example}[theorem]{Example}

\begin{document}
\title[Generalized point values in $\mathcal G(\mathbb Q_p^n)$]{On the characterization of $p$-adic Colombeau-Egorov generalized functions by their point values}
\author[E.\ Mayerhofer]{Eberhard Mayerhofer}
\address{Department of Mathematics, University of Vienna, Nordbergstrasse 15, 1090 Vienna, Austria}
\thanks{Work supported by FWF research grants P16742-N04 and Y237-N13}
\email{eberhard.mayerhofer@univie.ac.at}
\maketitle
\begin{abstract}
We show that contrary to recent papers by S.\ Albeverio, A.\ Yu.\ Khrennikov and V.\ Shelkovich, point values do not determine
elements of the so-called $p$-adic Colombeau-Egorov algebra $\mathcal G(\mathbb Q_p^n)$ uniquely. We further show in a more general way that for an Egorov algebra $\mathcal G(M,R)$ of generalized functions on a locally compact ultrametric space $(M,d)$ taking values
in a non-trivial ring, a point value characterization holds if and only if $(M,d)$ is discrete. Finally, following an idea due to M.\ Kunzinger and M.\ Oberguggenberger, a generalized point value characterization of $\mathcal G(M,R)$ is given.
Elements of $\mathcal G(\mathbb Q_p^n)$ are constructed which differ from the $p$-adic $\delta$-distribution considered as an element of $\mathcal G(\mathbb Q_p^n)$, yet coincide on point values with the latter. 
\end{abstract}
\section*{Introduction}
Algebras of generalized functions have found a growing number of applications in 
analysis, geometry and mathematical physics over the past two decades (cf.\ 
\cite{C, E, MObook, DIANA, K}, as well as \cite{DHPV} for a recent unified approach).
A distinguishing feature (compared to spaces of distributions in the sense of
Schwartz) of Colombeau- and Egorov type algebras is the availability
of a generalized point value characterization for elements of such spaces (see \cite{MO1},
resp.\ \cite{KS} for the manifold setting). Such a characterization may be
viewed as a nonstandard aspect of the theory: for uniquely determining
an element of a Colombeau- or Egorov algebra, its values on classical ('standard')
points do not suffice: there exist elements which vanish on each classical
point yet are nonzero in the quotient algebra underlying the respective
construction. A unique determination can only be attained by taking into
account values on generalized points, themselves given as equivalence 
classes of standard points. This characteristic feature is re-encountered
in practically all known variants of such algebras of generalized functions.

It therefore came as a surprise when in a series of papers (\cite{AKS, AKS2}) 
it was claimed that, contrary to the above general situation, in $p$-adic Colombeau-Egorov algebras 
a general point value characterization using only standard points was available.
This short note is dedicated to a thorough study of (generalized) point value characterizations
of $p$-adic Colombeau-Egorov algebras and to showing that in fact also in the $p$-adic
setting classical point values do not suffice to uniquely determine elements of such.

In the remainder of this section we recall some material from (\cite{AKS, AKS2}), using notation 
from \cite{DIANA}. 
Let $\mathbb N$ be the natural 
numbers starting with $n=1$. For a fixed prime $p$, let $\mathbb Q_p$ denote the field of rational $p$-adic numbers. 
Let $\mathcal D(\mathbb Q_p^n)$ denote the linear space of locally constant complex valued functions on 
$\mathbb Q_p^n$ ($n\geq 1$) with compact support. Let further $\mathcal P(\mathbb Q_p^n)
:=\mathcal D(\mathbb Q_p^n)^{\mathbb N}$. $\mathcal P(\mathbb Q_p^n)$ is endowed with an algebra-structure by defining addition and multiplication of sequences componentwise. Let $\mathcal N(\mathbb Q_p^n)$ be the subalgebra of elements $\{(f_k)_k\}\subset\mathcal P(\mathbb Q_p^n)$ such that for any compact set $K\subset \mathbb Q_p^n$ there exists an $N\in\mathbb N$ such that $\forall\; 
x\in K\;\forall\; k\geq N: f_k(x)=0$. This is an ideal in $\mathcal P(\mathbb Q_p^n)$. The quotient algebra 
$\mathcal G(\mathbb Q_p^n):=\mathcal P( \mathbb Q_p^n)/\mathcal N(\mathbb Q_p^n)$ is called the $p$-adic Colombeau-Egorov algebra. 
Finally, so called Colombeau-Egorov generalized numbers $\widetilde {\mathcal C}$ are introduced in the following way: 
Let $\bar {\mathbb C}$ be the one-point compactification of $\mathbb C\cup\{\infty\}$.\\ Factorizing $\mathcal A=\bar{\mathbb C}^{\mathbb N}$  by the ideal $\mathcal I:=\{u=(u_k)_k\in\mathcal A\mid \,\exists N\in\mathbb N\,\forall \;k\geq N: u_k=0\}$ yields then the ring $\widetilde{\mathcal C}$ of Colombeau-Egorov generalized numbers. We replace $\bar{\mathbb C}$ by $\mathbb C$ and construct similarly $\mathcal C$, the ring of generalized numbers: Clearly, $\bar{\mathbb C}$ is not needed in this context, since representatives of elements $f\in\mathcal G(\mathbb Q_p^n)$ merely take on values in $\mathbb C^{\mathbb N}$. 
Let $f=[(f_k)_k]\in\mathcal G(\mathbb Q_p^n)$. It is clear that for a fixed $x\in\mathbb Q_p^n$,  the {\it point value of f at x}, $[(f_k(x))_k]$ is a well defined element of
${\mathcal C}$, i.e., we may consider $f$ as a map
\begin{equation}\label{pv}
f: \;\mathbb Q_p^n\rightarrow \mathcal C:\;\;x\mapsto f(x):=(f_k(x))_k+\mathcal I.
\end{equation}
Note that the above constitutes a slight abuse of notation: 
The letter $f$ denotes both a generalized function (an element of $\mathcal G(\mathbb Q_p^n)$)
and a mapping on $\mathbb Q_p^n$. \\Finally, let $A$ be a set and let $R$ be a ring. For $B\subset A,\;\theta\in R$ we call the characteristic function of $B$ the map $\chi_{B,\theta}:\; A\rightarrow R$  which is identically $\theta$ on $B$ 
and which vanishes on $A\setminus B$. Furthermore, if $\theta=1\in R$ we simply write $\chi_B=\chi_{B,1}$. 
\section{Point values and a counterexample}
The following statement is proven in Theorem 4.4 of \cite{AKS}:
{\it Let $f\in\mathcal G(\mathbb Q_p^n)$, then:
\[
f=0 \;\;\mbox{in}\;\;\mathcal G(\mathbb Q_p^n)\Leftrightarrow\forall\; x\in\mathbb Q_p^n: f(x)=0\;\;\mbox{in}\;\; \mathcal C
\]
}
However, inspired by (\cite{zAR}, p.\ 218) we construct the following counterexample to this claim, 
which shows that point values cannot uniquely determine elements in $\mathcal G(\mathbb Q_p^n)$ uniquely. 
For the sake of simplicity we assume that $n=1$. 
\begin{example}\label{example}\rm
For any $l\in\mathbb N$, set
$$
B_l:=\{x\in\mathbb Z_p: \vert x-p^l\vert<\vert p^{2l}\vert\}\subset
\{x\in\mathbb Z_p: \vert x\vert=\vert p^l\vert\}\,. 
$$
For any $i \in\mathbb N$, we set $f_i:=\chi_{B_i}$. Clearly $B_i\cap B_j=\emptyset$ whenever $i\neq j$ and since $f_i\in\mathcal D(\mathbb Q_p)$ for all natural numbers $i$,
$(f_i)_i$ is a representative of some $f\in\mathcal G(\mathbb Q_p)$. Now, for any $\alpha\in \mathbb Q_p$, $f(\alpha)=0$ in $\mathcal C$, since either $\alpha\in B_i$ for some $i\in\mathbb N$ (which implies that $f_j(\alpha)=0\;\forall\; j>i$)   or $\alpha\in\mathbb Q_p\setminus \bigcup B_i$, where each $f_i$ ($i\in\mathbb N$) is identically zero. Consider now the sequence $(\beta_i)_{i\geq1}\in\mathbb Z_p^{\mathbb N}$, where $\beta_i=p^i\;\forall\;i\in\mathbb N$. It follows that $f_i(\beta_i)=1\;\forall\; i\in\mathbb N$. In particular, for $K=\mathbb Z_p$ or any clopen ball containing $0$, there is no representative $(g_j)_j$ of $f$ such that  for some $N>0$, $g_j=0\; \forall\; j\geq N$. Hence $f\neq 0$ in $\mathcal G(\mathbb Q_p)$ although all standard point values of $f$ vanish.
\end{example}
\begin{remark}
By means of the above example we may analyze the proof of Theorem 4.4 in \cite{AKS}. Let $f$ be the generalized function from \ref{example}. As a compact set choose $K:=B_{\leq p^{-2}}(0)=p^2\mathbb Z_p$. For the representative $(f_k)_k$ constructed in \ref{example} and $x=0$ we have $N(0)=1$, which in the notation of \cite{AKS} means that for any $k\geq 1=N(0)$, $f_k(0)=0$. Also, recall that $B_{\gamma}(a)$ is the dressed ball
$B_{\leq p^{\gamma}}(a)$. The ``parameter of constancy'' (\cite{AKS}, p.\ 6) of $f_1$ at $x=0$, which is the maximal $\gamma$ such that $f_1$ is identically zero on $B_{\gamma}(0)$, is $l_0(0)=-2$. Now, there exists a covering of $K$ consisting of a single set, namely $B_{l_0(0)}(0)$. Thus we may replace the application of the Heine-Borel Lemma in \cite{AKS} by our singleton-covering. But then the claim that (4.1) and (4.2) imply that for all $k\geq N(0)=1$ we have $f_k(0+x')=f_k(0)=0\;\forall\, x'\in K$ does not hold. This indeed follows from the definition of the sequence $(f_k)_k$ of locally constant functions from above, since for any $k\in\mathbb N$ we have $f_k(p^k)=1$. 
\end{remark}
\section{Egorov algebras on locally compact ultrametric spaces}
In this section we consider the problem of point value characterization in Egorov algebras in full generality: to this end we consider a
general locally compact ultrametric space $(M,d)$ instead of $\mathbb Q_p^n$, where $M$ need not have a field structure. Our aim is to show that even in such a general setting, the respective algebra cannot have a point value characterization, unless $M$ carries the discrete topology. Denote by $\mathcal E_d(M)$ the algebra of sequences of locally constant functions with compact support, taking values in a commutative ring $R\neq\{0\}$. Let $\mathcal N_d(M)$ be the set of negligible functions $\{(f_k)_k\}\subset\mathcal E_d(M)$
such that for any compact set $K\subset M$ there exists an $N\in\mathbb N$ such that $\forall \;x\in K\;\forall\; k\geq N: f_k(x)=0$. The subset $\mathcal N_d(M)$ is an ideal in $\mathcal E_d(M)$ and the quotient algebra $\mathcal G(M, R):=\mathcal E_d(M)/\mathcal N_d(M)$ is called the ultrametric Egorov algebra associated with $(M,d)$. Furthermore, the ring of generalized numbers is defined by $\mathcal R:=\mathbb R^{\mathbb N}/\sim$, where $\sim ~$ is the equivalence relation on $\mathbb R^{\mathbb N}$ given by
\[
u\sim v\;\mbox{in}\;\mathbb R^{\mathbb N}\Leftrightarrow \exists\,N\in\mathbb N\;\forall\;k\geq N: u_k-v_k=0.
\]
We call $\mathcal I(R):=\{w\in R^{\mathbb N}:w\sim 0\}$ the ideal of negligible sequences in $R$. Analogous to (\ref{pv}), for $f\in\mathcal G(M,R)$ evaluation on standard point values is introduced by means of the mapping:
\begin{equation}\label{pv1}
f: \;M\rightarrow \mathcal R:\;\;x\mapsto f(x):=(f_k(x))_k+\mathcal I(R).
\end{equation}
\begin{definition}
An ultrametric Egorov algebra $\mathcal G(M,R)$ is said to admit a standard point value characterization if for each $u\in\mathcal G(M,R)$ we have 
\[
u=0\Leftrightarrow \forall\;x\in M: u(x)=0\;\mbox{in}\;\mathcal R
\]
\end{definition}
Using this terminology, Example \ref{example} shows that $\mathcal G(\mathbb Q_p^n)$ does not admit a standard point value characterization. The main result of this section is:
\begin{theorem}
Let $(M,d)$ be a locally compact ultrametric space and let $R\neq \{0\}$. Then $\mathcal G(M, R)$ does not admit a standard point value characterization unless $(M,d)$ is discrete.
\end{theorem}
\begin{proof}
The result follows by generalizing the construction of Example \ref{example}. Assume $(M,d)$ is not discrete, then there exists a point $x\in M$ and a sequence $(x_n)_n$ of distinct points in $M$ converging to $x$. We may assume that $d(x,x_i)>d(x,x_j)$ whenever $i<j$. Define stripped balls $(B_n)_{n\geq 1}$ with centers $(x_n)_{n \geq 1} $ by $B_n:=\{y\in M\mid d(x_n,y)<\frac{d(x_n,x)}{2}\}$. Due to the ultrametric property ``the strongest one wins'' we have $B_n\subset \{z\mid d(x,z)=d(x_n,x)\}$, which further implies that for all $i\neq j,\,i,j\in\mathbb N$ the balls $B_i$, $B_j$ are disjoint sets in $M$. Since $R$ is a non-trivial ring, we may choose some $\theta\in R\setminus\{0\}$. Now we define a sequence $(f_k)_k$ of locally constant functions in the following way: For any $i\geq 1$ set $f_i=\chi_{B_i,\theta}$. Clearly, $f:=[(f_i)_i]\in\mathcal G(M, R)$, and similarly to Example \ref{example}, for any $\alpha\in M$, $f(\alpha)=0$  in $\mathcal R$. Nevertheless for the sequence $(x_n)_n$, which without loss of generality may be assumed to lie in a compact neighborhood of $x$, one has $f_i(x_i)=\theta\; \forall\; i\geq 1$ which implies that $f\neq 0$ in $\mathcal G(M,R)$.
\end{proof}
Recall that a discrete topological space $X$ has the following properties:
\begin{enumerate}
\item $X$ is locally compact.
\item Any compact set in $X$ contains finitely many points only.
\end{enumerate}
Therefore we know that for a set $D$ endowed with the discrete metric and for any commutative ring $R$, the respective ultrametric Egorov algebra $\mathcal G(D, R)$ admits a pointwise characterization. We therefore conclude:
\begin{corollary}
For a locally compact ultrametric space $(M,d)$ and a non-trivial ring $R$, the following statements are equivalent:
\begin{enumerate}
\item $\mathcal G(M, R)$ admits a standard point value characterization.
\item The topology of $(M,d)$ is discrete.
\end{enumerate}
\end{corollary}
\section{Generalized point values}
In this section we give an appropriate generalized point value characterization in the style of (\cite{DIANA}, pp.\ 37--43) of
$\mathcal G(M, R)$, where $M$ is endowed with a non-discrete ultrametric $d$ for which $M$ is locally compact, and $R\neq \{0\}$. First, we have to introduce a set $\widetilde M_c$ of compactly supported generalized points over $M$. Let $\mathcal E=M^{\mathbb N}$, the ring of sequences in $M$, and identify
two sequences, if for some index $N\in\mathbb N$ one has $d(x_n,y_n)=0\;\forall\; n\geq N$, i.e., $x_n=y_n\; \forall\; n\geq N$; we write $x\sim y$. We call $\widetilde M=\mathcal E/\sim$ the ring of generalized numbers. Finally, $\widetilde M_c$ is the subset
of such elements $x\in \widetilde M$ for which there exists a compact subset $K$  and some representative $(x_n)_n$ of $x$ such that for some $N>0$ we have $x_n\in K$ for all $n\geq N$. It follows that evaluating a function $u\in\mathcal G(M,R)$ at a compactly supported generalized point $x$ is possible, i.e., for representatives $(x_k)_k$, $(u_k)_k$ of $x$ resp.\  $u$, $[(u_k(x_k))_k]$ is a well defined element
of $\mathcal R$. 
\begin{proposition}
In $\mathcal G(M,R)$, there is a generalized point value characterization, i.e.,
\[
u=0 \;\mbox{in}\;\mathcal G(M,R)\;\Leftrightarrow\; \forall\; x\in \widetilde M_c: u(x)=0\;\; \mbox{in}\;\; \mathcal R
\]
\end{proposition}
\begin{proof}
The condition is obviously necessary. Conversely, let $u\in\mathcal G(M,R)$, $u\neq 0$. This means that there is a representative $(u_k)_k$ of $u$ and a compact set $K\subset\subset M$ such that $u_k$ does not vanish on $K$ for infinitely many $k\in\mathbb N$. In particular this means we have a sequence $(x_k)_k$ in $K$ such that for infinitely many $k\in\mathbb N$, $u_k(x_k)\neq 0$. Clearly this means that $u(x)\neq 0$ in $\mathcal R$ for the compactly supported generalized point defined by $x:=[(x_k)_k]$.
\end{proof}
\section{The imbedded $\delta$-distribution in $\mathcal G(\mathbb Q_p^n)$.}
In \cite{AKS},  Theorem 4.4 is illustrated by some examples, to highlight the advantage of a point value concept in $\mathcal G(\mathbb Q_p^n)$. In this section we discuss the $\delta$-distribution (Example 4.5 on p.\ 12 in \cite{AKS}) and construct a generalized function $f\in\mathcal G(\mathbb Q_p)$ different from $\delta$ which however coincides with $\delta$ on all standard points in $\mathbb Q_p$. We first imbed the $\delta$-distribution in $\mathcal G(\mathbb Q_p)$ as in \cite{AKS} (p.\ 9, Theorem 3.3) which yields $\iota(\delta)=(\delta_k)_k+\mathcal N_p(\mathbb Q_p)$, where $\delta_k(x):=p^k\Omega(p^k\vert x\vert_p)$ and
$\Omega$ is the bump function on $\mathbb R^+_0$ given by
\[
\Omega(t):=\begin{cases} 1,\;\;0\leq t\leq 1\\ 0,\;\;t>1\end{cases}
\]
Evaluation of $\iota(\delta)$ on standard points is shown in Example 4.5 of \cite{AKS}. With $\tilde c:=(p^k)_k+\mathcal I\in\mathcal C$
one has:
\[
\iota(\delta)(x)=\begin{cases}\widetilde c,\;\;x=0\\ 0,\;\;x\neq 0\end{cases}\qquad (x\in\mathbb Q_p)
\]
Let $\varphi: \mathbb N\rightarrow \mathbb Z$ be a monotonous function such that $\lim_{k\rightarrow\infty}\varphi(k)=\infty$, and
such that the cardinality of $U_{\varphi}:= \{k:\varphi(k)>k\}$ is infinite. Consider an element $f\in\mathcal G(\mathbb Q_p)$ given by $f:=(f_k)_k+\mathcal N(\mathbb Q_p)$ where for any 
$k\geq 1$, $f_k(x):=p^k\Omega(p^{\varphi(k)}\vert x\vert_p)$. Then the standard point values of $\iota(\delta)$ and $f$ coincide. Furthermore, they coincide on compactly supported generalized points $x\in\widetilde{\mathbb Q}_{p,c}$ with the property that for any representative $(x_k)_k$ of $x$ there exists an $N\in\mathbb N$ such that $\forall\;k\geq N:\;\vert x_k\vert_p>p^{-\min\{k,\varphi(k)\}}$, since in this case we have 
$\delta_k(x_k)=f_k(x_k)=0$. However, there are compactly supported generalized points
violating this condition which yield different generalized point values of $\iota(\delta)$ resp. $f$: consider the generalized point $x_0:=[(p^k)_k]\in\widetilde{\mathbb Q}_{p,c}$. Then $f(x_0)\neq \widetilde c$, since $\theta_k:=f_k(x_k)=0$ for any $k\in U_{\varphi}$ and thus $\theta_k=0$ for infinitely many $k\in\mathbb N$. But $\iota (\delta)(x_0)=\widetilde c$.

\end{document}